\documentclass[a4paper, 12pt]{amsart}

\usepackage{amsmath, amsfonts, amssymb, amsthm}

\newtheorem{thm}{Theorem}[section]
\newtheorem{lemma}[thm]{Lemma}

\providecommand{\aut}{\mathop{\rm Aut \,}\nolimits}

\begin{document}

\title{\bf Infinite primitive directed graphs}

\author{Simon Smith}



\subjclass{05 C 20, 20 B 15}

\date{\today}

\begin{abstract} A group $G$ of permutations of a set $\Omega$ is {\em primitive}
if it acts transitively on $\Omega$, and the only $G$-invariant
equivalence relations on $\Omega$ are the trivial and universal
relations. A graph $\Gamma$ is {\em primitive} if its automorphism
group acts primitively on its vertex set.

A graph $\Gamma$ has {\em connectivity one} if it is connected and
there exists a vertex $\alpha$ of $\Gamma$, such that the induced
graph $\Gamma \setminus \{\alpha\}$ is not connected. If $\Gamma$
has connectivity one, a {\em block} of $\Gamma$ is a connected
subgraph that is maximal subject to the condition that it does not
have connectivity one.

The primitive undirected graphs with connectivity one have been
fully classified by Jung and Watkins: the blocks of such graphs
are primitive, pairwise-isomorphic and have at least three
vertices. When one considers the general case of a directed
primitive graph with connectivity one, however, this result no
longer holds. In this paper we investigate the structure of these
directed graphs, and obtain a complete characterisation.

\end{abstract}

\maketitle

\section{Preliminaries}

Throughout this note, a {\it graph} will be a directed graph
without multiple edges or loops. A graph $\Gamma$ will be thought
of as a pair $(V\Gamma, E\Gamma)$, where $V\Gamma$ is the set of
{\it vertices} and $E\Gamma$ the set of {\it edges} of $\Gamma$.
The set $E\Gamma$ consists of ordered pairs of distinct elements
of $V\Gamma$. Two vertices $\alpha, \beta \in V\Gamma$ are {\it
adjacent} if either $(\alpha, \beta)$ or $(\beta, \alpha)$ lies in
$E\Gamma$. All paths will be undirected, unless otherwise stated.
A graph is {\em infinite} if its vertex set is infinite.

Two vertices are {\em connected} if there exists an undirected
path in $\Gamma$ between them, while a graph is {\it connected} if
any two vertices are connected. The {\em distance} between two
connected vertices $\alpha$ and $\beta$ in $\Gamma$ is the length
of the shortest path between them, and is denoted by
$d_\Gamma(\alpha, \beta)$.

We shall be deducing the properties of directed graphs from
characterisations of {\it undirected graphs}. This, of course,
requires us to carefully define the latter in a way that preserves
the natural relationship between the two. The term {\it undirected
graph} will henceforth refer to a directed graph with the property
that, whenever $(\alpha, \beta) \in E\Gamma$, we have $(\beta,
\alpha) \in E\Gamma$. In this case, it is sometimes convenient to
replace each pair of edges $(\alpha, \beta)$ and
$(\beta, \alpha)$ with the unordered pair $\{\alpha, \beta\}$.\\

Groups, and in particular groups of automorphisms, will play a
leading role in many of the arguments presented herein. Throughout
this work, $G$ will be a group of permutations of an infinite set
$\Omega$, where $\Omega$ will usually be the vertex set of some
infinite graph.

If $\alpha \in \Omega$ and $g \in G$, we denote the image of
$\alpha$ under $g$ by $\alpha^g$. Following this notation, all
permutations will act on the right. The set of images of $\alpha$
under all elements of $G$ is called an {\it orbit} of $G$, and is
denoted by $\alpha^G$. There is a natural action of $G$ on the
$n$-element subsets and $n$-tuples of $\Omega$ via $\{\alpha_1,
\ldots, \alpha_n\}^g := \{\alpha_1^g, \ldots, \alpha_n^g\}$ and
$(\alpha_1, \ldots, \alpha_n)^g := (\alpha_1^g, \ldots,
\alpha_n^g)$ respectively.

If $\alpha \in \Omega$, we denote the stabiliser of $\alpha$ in
$G$ by $G_\alpha$, and if $\Sigma \subseteq \Omega$ we denote the
setwise and pointwise stabilisers of $\Sigma$ in $G$ by
$G_{\{\Sigma\}}$ and $G_{(\Sigma)}$ respectively.

The group $G$ is {\em transitive} on $\Omega$ if $G$ has one orbit
on $\Omega$, namely $\Omega$ itself. A transitive group $G$ is
said to act {\em regularly} on $\Omega$ if $G_\alpha = 1$ for each
$\alpha \in \Omega$.

A {\em $G$-congruence on $\Omega$} is an equivalence relation
$\approx$ on $\Omega$ satisfying
\[\alpha \approx \beta \Leftrightarrow \alpha^g \approx \beta^g,\]
for all $\alpha, \beta \in \Omega$ and $g \in G$. A transitive
group $G$ is {\it primitive} on $\Omega$ if the only
$G$-congruencies admitted by $\Omega$ are the trivial and universal
equivalence relations. The following is well known.

\begin{thm} {\normalfont(\cite[Theorem 4.7]{neumann:ipg})} \label{thm:PrimAndMaxSubgps} If $G$
is a transitive group of permutations on $\Omega$, and $|\Omega| >
1$, then $G$ is primitive on $\Omega$ if and only if, for every
$\alpha \in \Omega$, the stabiliser $G_\alpha$ is a maximal
subgroup of $G$. \qed \end{thm}

A permutation $\sigma$ of $V\Gamma$ is an {\em automorphism} of
$\Gamma$ if it preserves the edge-structure of $\Gamma$; that is,
\[e \in E\Gamma \Leftrightarrow e^\sigma \in E\Gamma.\]
The set of all automorphisms of the graph $\Gamma$ form a group
called the {\em automorphism group of $\Gamma$}, denoted by $\aut
\Gamma$. A graph is {\em primitive} if $\aut \Gamma$ is primitive
on the set $V\Gamma$, and is {\em automorphism-regular} if $\aut
\Gamma$ acts regularly on $V\Gamma$.

The following theorem due to D.~G.~Higman gives a test for
primitivity.

\begin{thm} {\normalfont(\cite{higman})} \label{thm:higman} A
transitive group of permutations of a set $\Omega$ is
primitive if and only if every graph with vertex set $\Omega$ and
edge-set $(\alpha, \beta)^G$ is connected whenever $\alpha$ and
$\beta$ are distinct elements of $\Omega$. \qed \end{thm}

A useful consequence of this result is that every primitive graph must be connected.

The {\it connectivity} of an infinite connected graph $\Gamma$ is
the smallest possible size of a subset $W$ of $V\Gamma$ for which
the induced graph $\Gamma \setminus W$ is disconnected. A {\it
block} of $\Gamma$ is a connected subgraph that is maximal subject
to the condition it has connectivity strictly greater than one. If
$\Gamma$ has connectivity one, then the vertices $\alpha$ for
which $\Gamma \setminus \{\alpha\}$ is disconnected are called the
{\it cut vertices} of $\Gamma$.

A group acting on a graph $\Gamma$ is said to be {\it
vertex-transitive} or {\it edge-transitive} if it acts
transitively on the set of vertices or edges of $\Gamma$
respectively. Similarly, a graph admitting such a group will be
referred to as being {\em vertex-} or {\em edge-transitive}.

\section{Local structure}

Consider the following construction. Let $V_1$ be the set of cut
vertices of a connected graph $\Gamma$, and let $V_2$ be a set in
bijective correspondence with the set of blocks of $\Gamma$. We
let $T$ be a bipartite graph whose parts are $V_1$ and $V_2$. Two
vertices $\alpha \in V_1$ and $x \in V_2$ are adjacent in $T$ if
and only if $\alpha$ is contained in the block of $\Gamma$
corresponding to $x$; in this case, there are two edges in $T$
between the vertices $\alpha$ and $x$, one going in each
direction. Thus, $T$ can be considered to be an undirected graph.
In fact, this construction yields a tree, which is called the {\it
block-cut-vertex tree} of $\Gamma$. Note that if $\Gamma$ has
connectivity one and block-cut-vertex tree $T$, then any group $G$
acting on $\Gamma$ has a natural action on $T$.\\

Let $\Gamma$ be a primitive graph with connectivity one whose
blocks have at least three vertices, and suppose $G$ is a vertex-
and edge-transitive group of automorphisms of $\Gamma$. Since
$\Gamma$ is a vertex-transitive graph with connectivity one, every
vertex is a cut vertex. Fix some block $\Lambda$ of $\Gamma$, and
let $H$ be the subgroup of the automorphism group $\aut \Lambda$
induced by the setwise stabiliser $G_{\{\Lambda\}}$ of $V\Lambda$
in $G$. Let $T$ be the block-cut-vertex tree of $\Gamma$, and let
$x$ be the vertex of $T$ that corresponds to the block $\Lambda$.
Our aim in this section is to show $H$ is primitive but not
regular.

If $x_1$ and $x_2$ are distinct vertices of the tree $T$, we use
$C(T \setminus \{x_1\}, x_2)$ to denote the connected component of
$T \setminus \{x_1\}$ that contains the vertex $x_2$.

\begin{lemma} \label{lemma:block_stabiliser_is_primitive} If $G$ acts
primitively on the vertices of $\Gamma$, then the group $H$ acts
primitively on the vertices of $\Lambda$. \end{lemma}

\begin{proof} Fix $\alpha \in V\Lambda$ and suppose, for a
contradiction, the group $H$ does not act primitively on
$V\Lambda$. Then there is a vertex $\gamma \in V\Lambda \setminus
\{\alpha\}$ such that the graph $\Lambda^\prime := (V\Lambda,
(\alpha, \gamma)^H)$ is not connected. Let $\Gamma^\prime :=
(V\Gamma, (\alpha, \gamma)^G)$. We will show this graph cannot be
connected, and hence $G$ cannot be primitive.

Let $\{\Delta_i\}_{i \in I}$ be the set of connected components of
$\Lambda^\prime$ and let
\[\mathcal{C}_i := \bigcup_{\delta \in \Delta_i} C(T \setminus \{x\}, \delta) \cap V\Gamma.\]
Suppose $\delta_i \in \mathcal{C}_i$ and $\delta_j \in
\mathcal{C}_j$, with $i \not = j$. We claim $\delta_i$ and
$\delta_j$ are not adjacent in $\Gamma^\prime$. Indeed, since the
distance $d_T(\alpha, \gamma)$ between $\alpha$ and $\gamma$ in
$T$ is equal to $2$, if $\delta_i$ and $\delta_j$ are to be
adjacent in the edge-transitive graph $\Gamma^\prime$, it must be
the case that $d_T(\delta_i, \delta_j)=2$. If either $\delta_i$ or
$\delta_j$ is not adjacent to $x$ in $T$ then $d_T(\delta_i,
\delta_j) > 2$, so they cannot be adjacent in $\Gamma^\prime$. On
the other hand, if $\delta_i$ and $\delta_j$ are adjacent to $x$
in $T$, then they both lie in $V\Lambda = V\Lambda^\prime$, and
therefore $\delta_i \in \Delta_i$ and $\delta_j \in \Delta_j$. In
this case, if they are adjacent in $\Gamma^\prime$ then there
exists $g \in G$ such that either $(\delta_i, \delta_j)$ or
$(\delta_j, \delta_i)$ is equal to $(\alpha, \gamma)^g$. Such an
automorphism must fix $V\Lambda$ setwise, and therefore lies in
$G_{\{\Lambda\}}$. Thus, there exists an element $h \in H$ such
that either $(\delta_i, \delta_j)$ or $(\delta_j, \delta_i)$ is
equal to $(\alpha, \gamma)^h$, meaning that $\delta_i$ and
$\delta_j$ are adjacent in $\Lambda^\prime$; however, this
contradicts the fact that $\delta_i$ and $\delta_j$ are in
distinct components of $\Lambda^\prime$. Hence, $\delta_i$ and
$\delta_j$ are not adjacent in $\Gamma^\prime$.

Therefore, there can be no path in $\Gamma^\prime$ between a
vertex in $\mathcal{C}_i$ and a vertex in $\mathcal{C}_j$ whenever
$i \not = j$, and so the graph $\Gamma^\prime$ is not connected.
Whence, $G$ cannot act primitively on $V\Gamma$ by
Theorem~\ref{thm:higman}.
\end{proof}

Fix distinct vertices $\alpha, \beta \in V\Gamma$ and recall that
$\alpha$ and $\beta$ are also vertices of the block-cut-vertex
tree $T$.

A {\em geodesic} between two vertices is a shortest path between
them. In a tree, there is a unique geodesic between any two
vertices. Let $[\alpha, \beta]_T$ be the $T$-geodesic between
$\alpha$ and $\beta$, and let $(\alpha, \beta )_T$ be the
$T$-geodesic $[\alpha, \beta]_T$ excluding both $\alpha$ and
$\beta$. This notation extends obviously to $[\alpha, \beta)_T$
and $(\alpha, \beta]_T$.

Since $\alpha$ and $\beta$ are vertices of both $\Gamma$ and $T$,
the distance $d_T(\alpha, \beta)$ is even, so we may choose a
vertex $y \in (\alpha, \beta)_T$ that is distinct from $\alpha$
and $\beta$.

\begin{lemma} \label{technical_lemma:_components}
If $g \in G_\alpha$ does not fix $y \in VT$, and $\delta \not \in
C(T \setminus \{y\}, \alpha)$, then $\delta^g \not \in C(T
\setminus \{y\}, \beta)$. Similarly, if $g \in G_\beta$ does not
fix $y$ and $\delta \not \in C(T \setminus \{y\}, \beta)$ then
$\delta^g \not \in C(T \setminus \{y\}, \alpha)$. \end{lemma}

\begin{proof} If $\delta \not \in C(T \setminus \{y\}, \alpha)$ and
$\delta^g \in C(T \setminus \{y\}, \beta)$ then
$\delta, \delta^g \not \in C(T \setminus \{y\}, \alpha)$, so we
must have $g \in G_{\alpha, y}$. Similarly, if $\delta \not \in
C(T \setminus \{y\}, \beta)$ and $\delta^g \in C(T \setminus
\{y\}, \alpha)$ then $\delta, \delta^g \not \in C(T \setminus
\{y\}, \beta)$, so we must have $g \in G_{\beta, y}$. \end{proof}

\begin{lemma} \label{technical_lemma:_distance} If $g \in G_\alpha$ does
not fix the vertex $y$ and $\delta \not \in C(T \setminus
\{y\}, \alpha)$ then $d_T(y, \delta^g) > d_T(y, \delta)$.
Similarly, if $g \in G_\beta$ does not fix $y$ and $\delta \not
\in C(T \setminus \{y\}, \beta)$ then $d_T(y, \delta^g) > d_T(y,
\delta)$. \end{lemma}

\begin{proof} Let $y'$ be the vertex adjacent to $y$ in $[\alpha, y]_T$.
If $\delta \not \in C(T \setminus \{y\},
\alpha)$ then $y \in [\alpha, \delta]_T$. Since $g \in G_\alpha
\setminus G_y$, both $y$ and $y'$ lie on the geodesic $[\delta,
\delta^g]_T$, with $y' \in [\delta^g, y]_T$. Thus $d_T(\delta^g,
y) = d_T(\delta^g, y') + d_T(y', y)$. Now $d_T(\delta^g, y') \geq
d_T(\delta, y)+d_T(y, y')$, so $d_T(\delta^g, y) \geq d_T(\delta,
y)+1 > d_T(\delta, y)$. Interchanging $\alpha$ and $\beta$ in the
above argument completes the proof of this lemma. \end{proof}

Henceforth, if $H$ is a subgroup of $G$, then we will write $H
\leq G$; if we wish to exclude the possibility of $H = G$ we will
instead write $H < G$.

\begin{lemma} \label{technical_lemma:_action_of_G} Let $g_1, \ldots, g_n \in G_\alpha$
and $h_1, \ldots, h_n \in
G_\beta$, and suppose $G_{\alpha, y} = G_{\beta, y}$. If there
exists $\gamma \in V\Gamma$ such that $G_{\alpha, y} \leq
G_\gamma$ then, for some $m \leq n$, there exist $g_2^\prime,
\ldots, g_m^\prime \in G_\alpha \setminus G_y$ and $g_1^\prime \in
G_\alpha \setminus G_{y} \cup \{1\}$ together with $h_1^\prime,
\ldots, h_{m-1}^\prime \in G_\beta \setminus G_{y}$ and
$h_m^\prime \in G_\beta \setminus G_{y} \cup \{1\}$ such that
\[\gamma^{g_1^\prime h_1^\prime \ldots g_m^\prime h_m^\prime} = \gamma^{g_1 h_1 \ldots g_n h_n}.\] \end{lemma}

\begin{proof} The proof of this lemma will be an inductive argument. Suppose
there exists $\gamma \in VT$ such that $G_{\alpha, y} \leq G_{\gamma}$.

Let $n = 1$. When considering $h_1 \in G_\beta$ we have two cases:
either $h_1 \in G_{y}$ or $h_1 \in G_\beta \setminus G_{y}$. If
$h_1 \in G_{y}$ then $h_1 \in G_{\beta, y} = G_{\alpha, y}$, so
$g_1h_1 \in G_\alpha$. In this case, redefine $g_1:= g_1h_1$ and
set $h_1^\prime := 1$. Alternatively, if $h_1 \in G_\beta
\setminus G_{y}$ then set $h_1^\prime:=h_1$. Having found a
suitable $h_1^\prime$, we will now construct $g_1^\prime$ from the
(possibly redefined) element $g_1 \in G_\alpha$. We again have two
cases: either $g_1 \in G_{y}$ or $g_1 \in G_\alpha \setminus
G_{y}$. If $g_1 \in G_{y}$ then $g_1 \in G_{\alpha, y}$ and so
$g_1 \in G_{\gamma}$. In this case, we can choose $g_1^\prime :=
1$. Otherwise, if $g_1 \in G_\alpha \setminus G_{y}$, then choose
$g_1^\prime := g_1$. In choosing $g_1^\prime$ and $h_1^\prime$ in
this way we ensure
\[\gamma^{g_1h_1} = \gamma^{g_1^\prime h_1^\prime},\]
so the hypothesis holds when $n=1$.

Let $k$ be a positive integer, and suppose the hypothesis is true
for all integers $n \leq k$. Fix $g_1, \ldots, g_{k+1} \in
G_\alpha$ and $h_1, \ldots, h_{k+1} \in G_\beta$, and set
$\gamma^\prime := \gamma^{g_1 h_1 \ldots g_{k+1} h_{k+1}}$. We
will use induction to construct elements $g_2^\prime, \ldots,
g_m^\prime \in G_\alpha \setminus G_y$ and $g_1^\prime \in
G_\alpha \setminus G_{y} \cup \{1\}$ together with $h_1^\prime,
\ldots, h_{m-1}^\prime \in G_\beta \setminus G_{y}$ and
$h_m^\prime \in G_\beta \setminus G_{y} \cup \{1\}$ such that
\[\gamma^{g_1^\prime h_1^\prime \ldots g_m^\prime h_m^\prime} = \gamma^\prime,\]
where $m$ is some integer less than or equal to $k+1$.

We begin by considering $h_{k+1} \in G_\beta$. There are two
cases: either $h_{k+1} \in G_{y}$ or $h_{k+1} \in G_\beta
\setminus G_{y}$. If $h_{k+1} \in G_{y}$ then $h_{k+1} \in
G_{\beta, y} = G_{\alpha, y}$, so $g_{k+1}h_{k+1} \in G_\alpha$.
In this case, redefine $g_{k+1}:= g_{k+1}h_{k+1}$ and set
$h^\prime := 1$. If, on the other hand, $h_{k+1} \in G_\beta
\setminus G_{y}$, then set $h^\prime:=h_{k+1}$.

If we now consider the (possibly redefined) element $g_{k+1} \in
G_\alpha$, there are again two cases: either $g_{k+1} \in G_{y}$,
or $g_{k+1} \in G_\alpha \setminus G_{y}$. If $g_{k+1} \in G_{y}$
then $g_{k+1} \in G_{\alpha, y} = G_{\beta, y}$, so
$h_{k}g_{k+1}h^\prime \in G_\beta$. In this case, let $h'' :=
h_{k}g_{k+1}h^\prime$; then
\[\gamma^\prime = \gamma^{g_1 h_1 \ldots g_k h''},\]
so we can apply the induction hypothesis to $\gamma^{g_1 h_1
\ldots g_k h''}$ and we are done. If, on the other hand, $g_{k+1}
\in G_\alpha \setminus G_{y}$, then set $g^\prime := g_{k+1}$, and
observe
\[\gamma^\prime = \gamma^{g_1 h_1 \ldots g_k h_k g^\prime h^\prime}.\]
By the induction hypothesis, for some $m \leq k$ there exist
$g_2^\prime, \ldots, g_m^\prime \in G_\alpha \setminus G_{y}$ and
$g_1^\prime \in G_\alpha \setminus G_{y} \cup \{1\}$ together with
$h_1^\prime, \ldots, h_{m-1}^\prime \in G_\beta \setminus G_{y}$
and $h_m^\prime \in G_\beta \setminus G_{y} \cup \{1\}$ such that
\[\gamma^{g_1 h_1 \ldots g_k h_k} = \gamma^{g_1^\prime h_1^\prime \ldots g_m^\prime h_m^\prime}.\]
Set $g_{m+1}^\prime := g^\prime \in G_\alpha \setminus G_{y}$ and
$h_{m+1}^\prime := h^\prime \in G_\beta \setminus G_{y} \cup
\{1\}$. Then
\[\gamma^{g_1^\prime h_1^\prime \ldots g_{m+1}^\prime h_{m+1}^\prime} = \gamma^\prime,\]
so the hypothesis holds for $n = k+1$. \end{proof}

We are now in a position to present the main result of this section.

\begin{thm} \label{thm:Gya_eq_Gyb_not_prim} Let $G$ be a vertex-transitive group
of automorphisms of a connectivity-one graph $\Gamma$ whose
blocks have at least three vertices, and let $T$ be the
block-cut-vertex tree of $\Gamma$. If there exist distinct
vertices $\alpha, \beta \in V\Gamma$ such that, for some vertices
$\alpha', \beta' \in (\alpha, \beta)_T$,
\begin{enumerate}
\item \label{point1} $[\alpha, \alpha']_T \cap (\beta', \beta]_T = \emptyset$; and
\item \label{point2} $G_{\alpha, \alpha'} = G_{\beta, \beta'}$;
\end{enumerate}
then $G$ does not act primitively on $V\Gamma$.\end{thm}

\begin{proof} Suppose $G$ acts primitively on $V\Gamma$ and there exist distinct
vertices $\alpha, \beta \in V\Gamma$ and $\alpha', \beta'$ in the
$T$-geodesic $(\alpha, \beta)_T$, such that (\ref{point1}) and
(\ref{point2}) hold. We will show the group $\langle G_\alpha,
G_\beta \rangle$ generated by $G_\alpha$ and $G_\beta$ is not
transitive on $V\Gamma$. Then $G_\alpha < \langle G_\alpha,
G_\beta \rangle < G$, which, by applying
Theorem~\ref{thm:PrimAndMaxSubgps}, will contradict the assumption
that $G$ is primitive.

Choose $y \in [\alpha', \beta']_T$, and observe that by
(\ref{point2}) we have $G_{\alpha, y} = G_{\beta, y}$. Without
loss of generality, suppose $d_T(y, \alpha) \leq d_T(y, \beta)$.
As $G$ acts primitively on $V\Gamma$, the generated group $\langle
G_\alpha, G_\beta \rangle$ is not equal to $G_\alpha$, so we must
have $\langle G_\alpha, G_\beta \rangle = G$. Therefore the orbit
$\beta^{\langle G_\alpha, G_\beta \rangle}$ contains $\alpha$ and
there exist elements $g_1, \ldots, g_n \in G_\alpha$ and $h_1,
\ldots, h_n \in G_\beta$ such that $\alpha = \beta^{g_1h_1 \ldots
g_nh_n}$. By Lemma~\ref{technical_lemma:_action_of_G}, we can find
$g_2^\prime, \ldots, g_m^\prime \in G_\alpha \setminus G_{y}$ and
$g_1^\prime \in G_\alpha \setminus G_{y} \cup \{1\}$ together with
$h_1^\prime, \ldots, h_{m-1}^\prime \in G_\beta \setminus G_{y}$
and $h_m^\prime \in G_\beta \setminus G_{y} \cup \{1\}$ such that
\[\alpha = \beta^{g_1^\prime h_1^\prime \ldots g_m^\prime h_m^\prime}.\]
Suppose these automorphisms are chosen so that $m$ is minimal.

Now either $g_1^\prime \in G_\alpha \setminus G_{y}$ or
$g_1^\prime = 1$. If $g_1^\prime = 1$ then $\beta^{g_1'} = \beta$
and therefore $\beta^{g_1' h_1'} = \beta$. Thus $\beta^{g_2^\prime
h_2^\prime \ldots g_m^\prime h_m^\prime} = \alpha$, contradicting
the minimality of $m$. So we must have $g_1^\prime \in G_\alpha
\setminus G_{y}$. Since $\beta \not \in C(T \setminus \{y\},
\alpha)$, we may apply Lemma~\ref{technical_lemma:_components} and
Lemma~\ref{technical_lemma:_distance} to obtain $d_T(y,
\beta^{g_1^\prime}) > d_T(y, \beta)$ and $\beta^{g_1^\prime} \not
\in C(T \setminus \{y\}, \beta)$.

We now observe $h_1^\prime \not = 1$. Indeed, if $h_1^\prime = 1$
then $m=1$ and $\alpha = \beta^{g_1'}$; since $g_1' \in G_\alpha$
this is clearly not possible.

Thus, $h_1^\prime \in G_\beta \setminus G_{y}$ and
$\beta^{g_1^\prime} \not \in C(T \setminus \{y\}, \beta)$, and we
can again deduce from Lemma~\ref{technical_lemma:_components} and
Lemma~\ref{technical_lemma:_distance} that $d_T(y,
\beta^{g_1^\prime h_1^\prime}) > d_T(y, \beta^{g_1^\prime}) >
d_T(y, \beta)$, and $\beta^{g_1^\prime h_1^\prime} \not \in C(T
\setminus \{y\}, \alpha)$.

We may continue to apply Lemmas~\ref{technical_lemma:_components}
and \ref{technical_lemma:_distance} to obtain $\beta^{g_1^\prime
h_1^\prime \ldots g_m^\prime} \not \in C(T \setminus \{y\},
\beta)$ and $d_T(y, \beta^{g_1^\prime h_1^\prime \ldots
g_m^\prime}) > d_T(y, \beta)$. Now either $h_m^\prime \in G_\beta
\setminus G_{y}$ or $h_m^\prime = 1$. If $h_m^\prime = 1$, then
$\alpha = \beta^{g_1^\prime h_1^\prime \ldots g_m^\prime}$, so
$d_T(y, \alpha) = d_T(y, \beta^{g_1' h_1' \ldots g_m'})$ which is strictly
greater than $d_T(y,\beta)$. If $h_m^\prime \in G_\beta \setminus G_{y}$ then, by
Lemma~\ref{technical_lemma:_distance}, $d_T(y, \beta^{g_1^\prime
h_1^\prime \ldots g_m^\prime h_m^\prime}) > d_T(y, \beta)$; that
is, $d_T(y, \alpha) > d_T(y, \beta)$. Thus, in both cases $d_T(y,
\alpha) > d_T(y, \beta)$. This contradicts our assumption that
$d_T(y, \alpha) \leq d_T(y, \beta)$. Hence $\alpha \not \in
\beta^{\langle G_\alpha, G_\beta \rangle}$, and so $\langle
G_\alpha, G_\beta \rangle$ cannot act transitively on the set
$V\Gamma$.
\end{proof}

\begin{thm} \label{theorem:block_stabiliser_not_regular} Let $G$ be a
vertex-transitive group of automorphisms of a connectivity-one
graph $\Gamma$ whose blocks have at least three vertices. If $G$
acts primitively on $V\Gamma$ and $\Lambda$ is some block of
$\Gamma$ then $G_{\{\Lambda\}}$ is primitive and not regular on
$V\Lambda$. \end{thm}

\begin{proof} Suppose $\Lambda$ is a block of $\Gamma$ and $G_{\{\Lambda\}}$
acts primitively and regularly on $V\Lambda$. If $T$ is the
block-cut-vertex tree of $\Gamma$ then there exists a vertex $x
\in VT$ corresponding to the block $\Lambda$. Choose distinct
vertices $\alpha$ and $\beta$ in $V\Lambda$, and observe
$G_{\alpha, x} = G_{\alpha, \{\Lambda\}} \leq G_\beta$ and
$G_{\beta, x} = G_{\beta, \{\Lambda\}} \leq G_\alpha$.
Furthermore, $x \in (\alpha, \beta)_T$; hence $G$ cannot act
primitively on $V\Gamma$ by Theorem~\ref{thm:Gya_eq_Gyb_not_prim}.
\end{proof}

\section{Global structure}

In this section we will employ
Theorem~\ref{theorem:block_stabiliser_not_regular} to give a
complete characterisation of the primitive connectivity-one
directed graphs.

\begin{lemma} \label{lemma:H_acts_like_aut_Lambda} Suppose $\Gamma$ is a
vertex-transitive graph with connectivity one, whose blocks are
vertex-transitive, have at least three vertices and are pairwise
isomorphic. If $\Lambda$ is a block of $\Gamma$ and $H$ is the
subgroup of $\aut \Lambda$ induced by the action of $(\aut
\Gamma)_{\{\Lambda\}}$ on $\Lambda$, then $H = \aut \Lambda$.
\end{lemma}

\begin{proof} Let $T$ denote the block-cut-vertex tree of $\Gamma$, and let
$\Lambda$  be a block of $\Gamma$. We will show any automorphism
of the directed graph $\Lambda$ may be extended to an automorphism
of $\Gamma$.

We begin by asserting that if $\Lambda_1$ and $\Lambda_2$ are
blocks of $\Gamma$, and $\alpha_1$ and $\alpha_2$ are vertices in
$\Lambda_1$ and $\Lambda_2$ respectively, then there exists an
isomorphism $\rho : \Lambda_1 \rightarrow \Lambda_2$ such that
$\alpha_1^\rho = \alpha_2$. Indeed, by assumption, there exists an
isomorphism $\rho^\prime: \Lambda_1 \rightarrow \Lambda_2$. Define
$\alpha_1^\prime:= \alpha_1^{\rho^\prime}$. Since the block
$\Lambda_2$ is vertex-transitive, there exists an automorphism
$\tau$ of $\Lambda_2$ such that ${\alpha_1^{\prime}}^\tau =
\alpha_2$. Let $\rho:=\rho^\prime \tau$. Then $\rho:\Lambda_1
\rightarrow \Lambda_2$ is an isomorphism, with $\alpha_1^\rho =
\alpha_1^{\rho^\prime \tau} = {\alpha_1^\prime}^\tau = \alpha_2$.

Let $x$ be the vertex of $T$ that corresponds to $\Lambda$. For $k
\geq 0$, define $\Gamma_k$ to be the subgraph of $\Gamma$ induced
by the set $\{\alpha \in V\Gamma \mid d_T(x, \alpha) \leq 2k+1\}$.
We will show any automorphism $\sigma_k : \Gamma_k \rightarrow
\Gamma_k$ admits an extension $\sigma_{k+1}: \Gamma_{k+1}
\rightarrow \Gamma_{k+1}$. Whence, by induction, the lemma will
follow.

Fix $k \geq 0$ and let $\sigma_k : \Gamma_k \rightarrow \Gamma_k$
be an automorphism. Let $\{\alpha_i\}_{i \in I}$ be the set of
vertices in $V\Gamma_k \setminus V\Gamma_{k-1}$ (where
$V\Gamma_{-1} := \emptyset$). Each vertex $\alpha_i$ belongs to a
unique block $\Lambda_i$ of $\Gamma_k$, and, if $k \geq 1$, the
block $\Lambda_i$ possesses exactly one vertex in $\Gamma_{k-1}$.
Since $\Gamma$ is vertex transitive, any two vertices lie in the
same number of blocks of $\Gamma$, so let $\{\Lambda_{i, j}\}_{j
\in J}$ be the set of blocks of $\Gamma$ that contain $\alpha_i$
and are distinct from $\Lambda_i$. Each block $\Lambda_{i, j}$ is
wholly contained in $\Gamma_{k+1}$ and has exactly one vertex in
$\Gamma_k$, namely $\alpha_i$. If $i \in I$, set
$\alpha_i^\prime:=\alpha_i^{\sigma_k}$ and
$\Lambda_i^\prime:=\Lambda_i^{\sigma_k}$. Then $\Lambda_i^\prime =
\Lambda_{i^\prime}$ for some $i^\prime \in I$. For all $j \in J$
there exists an isomorphism $\rho_{i, j} : \Lambda_{i, j}
\rightarrow \Lambda_{i^\prime, j}$ such that $\alpha_i^{\rho_{i,
j}} = \alpha_i^\prime$. Thus, we may define a mapping
$\sigma_{k+1} : \Gamma_{k+1} \rightarrow \Gamma_{k+1}$ with
\[ \beta^{\sigma_{k+1}} :=
    \begin{cases}
    \beta^{\sigma_k} & \text{if $\beta \in V\Gamma_k$;} \\
    \beta^{\rho_{i, j}} & \text{if $\beta \in V\Lambda_{i, j}$.}
    \end{cases}
\]
This is clearly a well-defined automorphism of $\Gamma_{k+1}$.
\end{proof}

The primitive undirected graphs with connectivity one have the
following complete characterisation.

\begin{thm} {\normalfont(\cite[Theorem 4.2]{jung:watkins})}
\label{thm:undirected_connectivity_one_graphs} If $\Gamma$ is a
vertex-transitive undirected graph with connectivity one, then it
is primitive if and only if the blocks of $\Gamma$ are primitive,
pairwise isomorphic and each has at least three vertices. \qed
\end{thm}

This useful result seems to suggest a similar characterisation may
be possible for directed primitive graphs with connectivity one.
This is indeed the case.

\begin{thm} \label{thm:directed_connectivity_one_graphs} If $\Gamma$ is a
vertex-transitive directed graph with connectivity one, then it is
primitive if and only if the blocks of $\Gamma$ are primitive but
not automorphism-regular, pairwise isomorphic and each has at
least three vertices. \end{thm}

\begin{proof} Let $\Gamma$ be a directed vertex-transitive graph with
connectivity one. Suppose the blocks of $\Gamma$ are primitive but
not automorphism-regular, pairwise isomorphic and each has at
least three vertices. Let $\approx$ be a non-trivial $\aut
\Gamma$-congruence on $V\Gamma$. We will show this relation must
be universal, and thus that $\Gamma$ is a primitive graph. Since
the relation is non-trivial, there exist distinct vertices
$\alpha, \beta \in V\Gamma$ such that $\alpha \approx \beta$. Let
$T$ be the block-cut-vertex tree of $\Gamma$, let $\gamma \in
V\Gamma$ be the vertex in the geodesic $[\alpha, \beta]_T$ such
that $d_T(\beta, \gamma) = 2$, and let $\Lambda$ be the block of
$\Gamma$ containing $\beta$ and $\gamma$. By
Lemma~\ref{lemma:H_acts_like_aut_Lambda}, the group $(\aut
\Gamma)_{\{\Lambda\}}$ acts primitively but not regularly on
$V\Lambda$. Thus, there exists an automorphism $g \in (\aut
\Gamma)_{\gamma, \{\Lambda\}}$ that does not fix $\beta$. We are
considering the full automorphism group of the connectivity-one
graph $\Gamma$, so there must therefore exist an element $g^\prime
\in (\aut \Gamma)_{\alpha, \gamma, \{\Lambda\}}$ that does not fix
$\beta$. Thus, $\beta$ and $\beta^{g^\prime}$ are distinct
vertices in $\Lambda$. Now $\alpha \approx \beta$, so $\alpha
\approx \beta^{g\prime}$, and therefore $\beta \approx
\beta^{g^\prime}$. Since $(\aut \Gamma)_{\{\Lambda\}}$ is
primitive on $V\Lambda$ and $\approx$ induces a non-trivial $(\aut
\Gamma)_{\{\Lambda\}}$-congruence on $V\Lambda$, this relation
must be universal in $\Lambda$. By assumption, $\aut \Gamma$ acts
transitively on the blocks of $\Gamma$, so if two vertices lie in
the same block then they must lie in the same congruence class.
Thus, if $\gamma$ is any vertex of $\Gamma$, and $\alpha x_1
\alpha_1 x_2 \ldots x_n \gamma$ is the geodesic in $T$ between
$\alpha$ and $\gamma$, then $\alpha$ and $\alpha_1$ lie in a
common block, so $\alpha \approx \alpha_1$. Similarly, $\alpha_1
\approx \alpha_2$ and $\alpha_2 \approx \alpha_3$, so $\alpha
\approx \alpha_2$ and $\alpha \approx \alpha_3$. Continuing in
this way we eventually obtain $\alpha \approx \gamma$. Hence, this
congruence relation is universal on $V\Gamma$.

Conversely, suppose the group $\aut \Gamma$ acts primitively on
$V\Gamma$. Since $\Gamma$ is a directed primitive graph with
connectivity one, we can obtain an undirected graph
$\Gamma^\prime$ with vertex set $V\Gamma$ and edge set
$\{\{\alpha, \beta\} \mid (\alpha, \beta) \in E\Gamma\}$. Two
vertices are adjacent in $\Gamma$ if and only if they are adjacent
in $\Gamma^\prime$. As $\aut \Gamma$ is primitive on $V\Gamma$ and
$\aut \Gamma \leq \aut \Gamma^\prime$, it follows that $\aut
\Gamma^\prime$ must be primitive on $V\Gamma$, and hence
$\Gamma^\prime$ is a primitive undirected graph. Since $\Gamma$
has connectivity one, the same is true of $\Gamma^\prime$, so we
may apply Theorem~\ref{thm:undirected_connectivity_one_graphs} to
deduce the blocks of $\Gamma^\prime$ are primitive, pairwise
isomorphic and each has at least three vertices. Now, given a
block $\Lambda$ of $\Gamma$, there is a block $\Lambda^\prime$ of
$\Gamma^\prime$ such that $V\Lambda = V\Lambda^\prime$. Therefore,
the blocks of $\Gamma$ have at least three vertices, and are
primitive but not automorphism-regular by
Theorem~\ref{theorem:block_stabiliser_not_regular}.

It remains to show they are pairwise isomorphic. Fix some block
$\Lambda$ of $\Gamma$ and an edge $(\alpha, \beta) \in E\Lambda$.
Let $\Gamma_1$ be the graph $(V\Gamma, (\alpha, \beta)^{\aut
\Gamma})$. As $\aut \Gamma$ is primitive, this graph is a
connected subgraph of $\Gamma$. Thus, every block of $\Gamma$ must
contain an edge in $E\Gamma_1$. Furthermore, if $\Lambda^\prime$
is a block of $\Gamma$, then any automorphism of $\Gamma$ mapping
the edge $(\alpha, \beta)$ to an edge in $\Lambda^\prime$ must map
$\Lambda$ to $\Lambda^\prime$. Since $\Gamma_1$ is
edge-transitive, the blocks of $\Gamma$ must be pairwise
isomorphic. \end{proof}

This paper forms part of the author's DPhil thesis, completed
under the supervision of Peter Neumann at the University of
Oxford. The author would like to thank Dr Neumann for his tireless
enthusiasm and helpful suggestions. The author would also like to
thank the EPSRC for generously funding this research.


\vspace{3mm}

{\em Email address}: simonmarksmith@gmail.com

\end{document}